\documentclass[12pt,a4paper]{article}
\usepackage{amsmath,amsthm,amsfonts,amssymb,alltt}
\usepackage{epsfig}
\usepackage{color}
\usepackage{lscape}

\newtheorem{theorem}{Theorem}
\newtheorem{lemma}[theorem]{Lemma}

\theoremstyle{definition}

\theoremstyle{remark}
\newtheorem{remark}[theorem]{Remark}

\numberwithin{equation}{section}

\renewcommand{\Re}{\mathrm{Re}}
\renewcommand{\Im}{\mathrm{Im}}

\begin{document}
\sloppy

\begin{center}
{\bf Riemann boundary value problem with piecewise constant
matrix. \\
Part I. General algorithm}
\end{center}

 \vspace{5mm}

\begin{center}

Vladimir V. Mityushev

 \vspace*{5mm}
\begin{footnotesize}

Department of Computer Sciences, Pedagogical University,
\\
ul. Podchorazych 2, 30-084 Krakow, Poland

\vspace*{1mm}
 mityu@up.krakow.pl

\end{footnotesize}
\end{center}

\begin{abstract}
The vector-matrix  Riemann boundary value problem for the unit disk with piecewise constant matrix is constructively solved  by a method of functional equations. By functional equations we mean iterative functional equations  with shifts involving compositions of unknown functions analytic in mutually disjoint disks. The functional equations are written as an infinite linear algebraic system on the coefficients of the corresponding Taylor series. The compactness of the shift operators implies justification of the truncation method for this infinite system. The unknown functions and partial indices can be calculated by truncated systems. 
\end{abstract}

\section{Introduction}
Solution to the vector-matrix  Riemann boundary value problem (${\mathbb C}$--linear conjugation problem in terminology \cite{MiR}) and estimation of its partial indices is a fundamental problem of complex analysis. Though the general theory concerning solvability and existence of partial indices can be considered as completed \cite{vek} various special problems have been constructively solving in applications, see papers beginning from \cite{Khr71}. One can find a review and extensive references devoted to constructive solution to the problem in \cite{RogMit}, \cite{Antipov}.

The vector-matrix  Riemann boundary value problem for the unit disk with piecewise constant matrix attracts special attention as a model fundamental problem. It follows froim the Riemann-Hilbert monodromy problem (Hilbert 21st problem) where the partial indices play the crucial role \cite{Gia}-\cite{Gia3}. Exact solutions to the problem in special cases were discussed in \cite{ErhSpe02}, \cite{ErhSp01}, \cite{Khv94}, \cite{MeiSpe89}, \cite{RogMit}, \cite{ZveKhv85}.

The aim of this paper is a construction solution to this problem for an arbitrary matrix dimension and an arbitrary finite number of discontinuity points. The main idea of the method is the reduction of the problem to a system of functional equations in complex domain. By functional equations we mean iterative functional equations involving compositions of unknown functions with shifts without integral terms discussed in \cite{MiR}. The functional equations for analytic functions are constructed for a circular multiply connected domain. The compactness of the corresponding operators is established. Using the method of series we can find the unknown functions analytic in disks by their Taylor series. The coefficients of the Taylor series satisfy an infinite linear algebraic system. The compactness implies justification of the truncation method for this infinite system. The partial indices can be calculated by the truncated systems. 

In this first part, a general constructive procedure to determine the partial indices and to construct an approximate formula for the canonical matrix is described. The second part will  concern numerical implementation. 

Let $n$ points $\zeta_k$ divide the unit circle $|\zeta|=1$ onto
$n$ arcs $L_k=(\zeta_{k-1},\zeta_k)$. It is convenient to assume
that $k$ takes the values $1,2,\ldots,n$ which belong to the group
${\mathbb Z}_n$ (integers modulo $n$ where $n \equiv 0$). Let ${\mathbb D}^+$ be the unit disk $|\zeta|<1$ and ${\mathbb D}^-$ be its exterior, i.e.,
$|\zeta|>1$. We introduce a topology on the extended complex plane
$\widehat{{\mathbb C}}={\mathbb C} \cup \{\infty\}$ in such a way
that the closures of ${\mathbb D}^{\pm}$ do not intersect.

Let a constant matrix $G_k$ of dimension $N\times N$ ($N>1$) be given on each arc $L_k$. The following restrictions are imposed on the matrices $G_k$  ($k=1,2,\ldots,n$):\newline

{\bf Restriction 1.} Each matrix $G_k$ is invertible. \newline

{\bf Restriction 2.}
Each matrix $G_k^{-1}G_{k-1}$ has different eigenvalues.\newline

These restrictions reduce the number of different cases which should be separately discussed following the same lines for complete solution to the problem.

By the order at infinity of a scalar function meromorphic in
$|\zeta|>1$ we call the order of pole if it has a pole at
$\infty$, and the order of zero with the sign minus, if it is
analytic at $\infty$. In particular, the order is equal to zero,
if the function has a non-zero value at $\infty$. By the order of
a vector function at infinity we call the maximum of the orders of
its components.

The Riemann (${\mathbb C}$--linear conjugation) problem in the
theory of boundary value problems is stated as follows. To find a
vector function
$\Phi(\zeta)=(\Phi_1(\zeta),\Phi_2(\zeta),\ldots,\Phi_N(\zeta))^{\intercal}$
which is analytic in ${\mathbb D}^+ \cup {\mathbb D}^-$,
continuous in the closures of ${\mathbb D}^+$ and ${\mathbb D}^-$
except the points $V=\{\zeta_1,\zeta_2,\ldots, \zeta_n\}$ where it
is almost bounded (see the definition in the next section), with the following conjugation conditions
\begin{equation}
\label{2.1} \Phi^+(\zeta)=G_k \Phi^-(\zeta), \quad \zeta\in L_k, \quad
k=1,2,\ldots,n.
\end{equation}
Moreover, $\Phi(\zeta)$ has a finite order at infinity.
In relation (\ref{2.1}), $\Phi^+(\zeta)$ denotes the
boundary limit value $\lim_{\varsigma \to \zeta}\Phi(\varsigma)$,
where $|\varsigma|<1$ and $\Phi^-(\zeta) = \lim_{\varsigma \to \zeta}\Phi(\varsigma)$ where $|\varsigma|>1$. It is consistent with the positive
orientation of the unit circle. Vector functions having the
boundary values $\Phi^+(\zeta)$ and $\Phi^-(\zeta)$ from the left and the right sides of $|\zeta|=1$ can be also considered as two separate functions analytic in $|\varsigma|<1$ and in $|\varsigma|>1$, respectively.

The following general designations \cite{gakh, mus} are used. Let $D$ be a domain on the complex plane $\widehat{{\mathbb C}}$ and $U$ be a finite set
of the boundary points on $\partial D$. Introduce the class of
vector functions $\mathfrak{H}(D,U)$ analytic in $D$, continuous
in $D\setminus U$ and almost bounded at $U$. If the domain $D=\dot{D} \cup
\{\infty\}$ is unbounded, we distinguish the class of vector
functions $\mathfrak{H}(\dot{D},U)$ having a finite order at
infinity. Moreover, we distinguish one more class
$\mathfrak{H}_0(D,U)$ of vector functions which vanish at
infinity. Hence, the problem (\ref{2.1}) has been stated above in
the class $\mathfrak{H}({\mathbb{D}}^+\cup \dot{\mathbb{D}}^-,V)$.
The class $\mathfrak{H}(D,U)$ is introduced for $N$--dimensional vector functions. The same designation  are used for $2N$--dimensional and scalar functions.

\section{Classic theory} \label{sec2}

\subsection{Scalar problem} \label{sec2a}
\label{sec:sc}
In the present section, the well--known scalar Riemann
($\mathbb{C}$--linear conjugation) problem with discontinuous
coefficients are summarized due to F. D. Gakhov \cite{gakh} and N. I.
Muskhelishvili \cite{mus}. We consider only the case of piecewise
constant coefficient. Despite the results are known, they are
presented in non--traditional form which is needed in further
investigations. 

A scalar function $f(\zeta)$ analytic in $D$, continuous in its closure
except a point $\zeta_0 \in \partial D$ is called the almost
bounded at $\zeta=\zeta_0$ if $\lim_{\zeta \to \zeta_0}  |\zeta-\zeta_0|^{\varepsilon} f(\zeta) =0$ for arbitrary positive $\varepsilon$. We introduce the class of functions
$\mathfrak{G}(D,U)$ which differs from the class
$\mathfrak{H}(D,U)$ by the almost boundness instead of the
boundness at $U$.

Let $\nu_k$ ($k=1,2,\ldots,n$) be nonzero numbers, $b(\zeta)$ and
$b^{\prime}(\zeta)$ \footnote{Here, the prime is not a
derivative.} be functions H\"{o}lder continuous on $|\zeta|=1$
except a finite set of points $V$ where they have one--sided limits. The scalar Riemann
problem is stated as follows. To find $\Phi \in
\mathfrak{G}_0(\mathbb{D}^+\cup \mathbb{D}^-,V)$ with the
$\mathbb{C}$--conjugation condition
\begin{equation}
\label{g1} \Phi^+(\zeta)=\nu_k \Phi^-(\zeta)+b(\zeta), \; \zeta
\in L_k, \; k=1,2,\ldots,n.
\end{equation}
The problem (\ref{g1}) will be discussed simultaneously with the problem
\begin{equation}
\label{g1b} \Omega^+(\zeta)=\nu_k^{-1}
\Omega^-(\zeta)+b^{\prime}(\zeta), \; \zeta \in L_k, \;
k=1,2,\ldots,n.
\end{equation}

We begin the study from two homogeneous problems in
$\mathfrak{H}(\mathbb{D}^+\cup \dot{\mathbb{D}}^-,V)$
\begin{equation}
\label{g2} \Phi^+(\zeta)=\nu_k \Phi^-(\zeta), \; \zeta \in L_k, \;
k=1,2,\ldots,n,
\end{equation}
\begin{equation}
\label{g2b} \Omega^+(\zeta)=\nu_k^{-1} \Omega^-(\zeta), \; \zeta
\in L_k, \; k=1,2,\ldots,n.
\end{equation}
Let  $\theta_k\in (-\pi,\pi]$ be the argument of $\nu_{k-1}\nu_k^{-1}$, i.e., for a positive $\rho_k$,
\begin{equation}
\label{g4nu}
\nu_{k-1}\nu_k^{-1}=\rho_k e^{i\theta_k}.
\end{equation}
Introduce the logarithm
\begin{equation}
\label{g4a}
\gamma_k=\frac1{2\pi i} \ln \nu_{k-1}\nu_k^{-1}=
\frac{\theta_k}{2\pi}-\chi_k -\frac{i}{2 \pi}\ln \rho_k,
\end{equation}
where the integer $\chi_k$ is chosen in such a way that
\begin{equation}
\label{g4}
0 \leq \frac{\theta_k}{2\pi}-\chi_k<1.
\end{equation}
Along similar lines $\chi_k^{\prime}$ is chosen by the inequality
\begin{equation}
\label{g5} 0 \leq -\frac{\theta_k}{2\pi}-\chi_k^{\prime}<1.
\end{equation}
Note that $-\theta_k$ is the argument of $\nu_k\nu_{k-1}^{-1}$.
Let $\chi = \sum_{k=1}^n \chi_k$ and $\chi^{\prime} = \sum_{k=1}^n \chi_k^{\prime}$ be the indices of the problems (\ref{g2}) and
(\ref{g2b}), respectively.

Divide the set $V$ onto two disjoint subsets
$\mathcal{V}=\{\zeta_k \in V: \theta_k \neq 0\}$ and
$\mathcal{W}=\{\zeta_k \in V: \theta_k =0\}$ ($V=\mathcal{V} \cup
\mathcal{W}$). It follows from (\ref{g4})--(\ref{g5}) that
$\chi_k+\chi_k^{\prime}=-1$ for $\zeta_k \in \mathcal{V}$ and
$\chi_k+\chi_k^{\prime}=0$ for $\zeta_k \in \mathcal{W}$. Let $\#
(\mathcal{V})$ be the number of points in $
\mathcal{V}$, then
\begin{equation}
\label{g6} -(\chi_k+\chi_k^{\prime})=\# (\mathcal{V}).
\end{equation}

According to \cite{gakh} (p.437) the non--homogeneous problems
(\ref{g1}) and (\ref{g1b}) are solvable in
$\mathfrak{G}_0(\mathbb{D}^+\cup \mathbb{D}^-,V)$ if and only if
$\# (\mathcal{V})$ linearly independent conditions on the right
parts are fulfilled
\begin{equation}
\label{g7} \int_{|\tau|=1} [X_1^+(\tau)]^{-1} \prod_{k=1}^n
(\tau-\zeta_k)^{-\gamma_k} b(\tau) \tau^{j-1}d\tau=0, \;
j=1,2,\ldots,-\chi,
\end{equation}
\begin{equation}
\label{g7b} \int_{|\tau|=1} [X_1^-(\tau)]^{-1} \prod_{k=1}^n
(\tau-\zeta_k)^{\gamma_k} b^{\prime}(\tau) \tau^{j-1}d\tau=0, \;
j=1,2,\ldots,-\chi^{\prime},
\end{equation}
where
\begin{equation}
\label{g7c} X_1^+(\zeta)=\exp \Gamma^+(\zeta), \quad X_1^-(\zeta)=
z^{-\chi} \exp \Gamma^-(\zeta),
\end{equation}

\begin{equation}
\label{g8} \Gamma(\zeta)= \frac1{2\pi i} \sum_{k=1}^n \int_{|\tau|=1}
\ln\left[\nu_k \tau^{-\frac1{2\pi}  \sum_{m=1}^n (\theta_m-i\ln
\rho_m)}\right]\frac{d\tau}{\tau-\zeta}, \quad \zeta \in \mathbb{D}^+\cup \mathbb{D}^-.
\end{equation}
If (\ref{g7}) is fulfilled, the unique solution of (\ref{g1}) has
the form
\begin{equation}
\label{g9} \Phi^+(\zeta) = \prod_{k=1}^n
(\zeta-\zeta_k)^{\gamma_k} X_1^+(\zeta) \Psi_1^+(\zeta), \quad
\Phi^-(\zeta) = \prod_{k=1}^n
\left(\frac{\zeta-\zeta_k}{\zeta}\right)^{\gamma_k} X_1^-(\zeta)
\Psi_1^-(\zeta),
\end{equation}
where
\begin{equation}
\label{g9a} \Psi_1^+(\zeta) = \frac1{2\pi i} \int_{|\tau|=1}
\prod_{k=1}^n (\tau-\zeta_k)^{-\gamma_k} [X_1^+(\tau)]^{-1}
b(\tau)\frac{d\tau}{\tau-\zeta}.
\end{equation}
The solution of (\ref{g1b}) has analogous form.

It follows from (\ref{g7c})--(\ref{g9a}) that $\lim_{\zeta \to
\zeta_k} \Phi^{\pm} (\zeta) =0$ for $\zeta_k \in \mathcal{V}$. If
$\zeta_k \in \mathcal{W}$, then $\Re \gamma_k =0$ and
$|\zeta-\zeta_k|^{\gamma_k}$ is bounded as $\zeta$ tends to
$\zeta_k$. Moreover, $\Psi_1(\zeta)$ has a logarithmic singularity
at $\zeta=\zeta_k$, since $|\zeta-\zeta_k|^{-\gamma_k}$ is bounded
at $\zeta=\zeta_k$ and $b(\zeta)$ is discontinuous at
$\zeta=\zeta_k$. Similar asymptotic formulae are valid for
$\Omega^{\pm} (\zeta)$. Thus, we have
\begin{lemma}
\label{LG} The non--homogeneous problems (\ref{g1}) and
(\ref{g1b}) in $\mathfrak{H}_0(\mathbb{D}^+\cup \mathbb{D}^-,V)$
are solvable if and only if

i) $\# (\mathcal{V})$ linearly independent conditions (\ref{g7})
and (\ref{g7b}) are fulfilled,

ii) $b(\zeta)$ and $b^{\prime}(\zeta)$ are continuous at the
points of $\mathcal{W}$.
\end{lemma}

{\bf Remark}. 
The condition ii) can be written in the form $b(\zeta_k-0) =
b(\zeta_k+0)$, $k=1,2,\ldots, \#(\mathcal{W})$. This implies that
the problems (\ref{g1}) and (\ref{g1b}) in
$\mathfrak{H}_0(\mathbb{D}^+\cup \mathbb{D}^-,V)$ are solvable if
and only if $n=\#(\mathcal{V})+\#(\mathcal{W})$ conditions are
fulfilled.

\subsection{N. P. Vekua's theory} \label{sec2b}
The exhaustive theory of solvability of the vector--matrix problem
(\ref{2.1}) is presented in the book \cite{vek}. Though
\cite{vek} does not contain general exact formulas for
$\Phi(\zeta)$, Vekua's theory is useful in our constructive
method. In the present section, some results from \cite{vek}
adopted to our investigations are shortly outlined.

Let the invertible matrix function $G(\zeta)$ be H\"{o}lder
continuous on the unit circle except the set
$V=\{\zeta_1,\zeta_2,\ldots, \zeta_n\}$ where $G(\zeta)$ has
one--sided limits. Let $b(\zeta)$ be H\"{o}lder continuous on the
unit circle except $V$ where it has one--sided limits. Consider
the non--homogeneous Riemann (${\mathbb C}$--linear conjugation
problem) in the class $\mathfrak{H}_0({\mathbb{D}}^+\cup
{\mathbb{D}}^-,V)$
\begin{equation}
\label{v1} \Phi^+(\zeta)=G(\zeta) \Phi^-(\zeta)+b(\zeta), \;
|\zeta|=1,
\end{equation}
and the homogeneous problem in the class
$\mathfrak{H}({\mathbb{D}}^+\cup \dot{\mathbb{D}}^-,V)$
\begin{equation}
\label{v2} \Phi^+(\zeta)=G(\zeta) \Phi^-(\zeta), \; |\zeta|=1.
\end{equation}
The order at infinity, $M$, of a solution of (\ref{v2}) is
considered here as a control parameter. If for each
$k=1,2,\ldots,n$ the matrix $G(\zeta_k+0)^{-1}G(\zeta_k-0)$ has
different eigenvalues, the problem (\ref{v2}) is solvable for
sufficiently large $M$ (see Section 18, \cite{vek}).

By the linear independence of solutions of (\ref{v2}) we mean the
linear independence with polynomial coefficients. Since for
sufficiently large $M$ the problem (\ref{v2}) is solvable, hence
among its solutions there exist some with the lowest possible $M$.
Denote this order by $-\kappa_1$ and by $X_1(\zeta)$ a solution of
(\ref{v2}) of order $-\kappa_1$. The next solution $X_2(\zeta)$
has the lowest possible order $-\kappa_2$ from all solution of
(\ref{v2}) linearly independent with $X_1(\zeta)$. From the
remaining solutions, consider all those are linearly independent
with $X_1(\zeta)$, $X_2(\zeta)$. Among these solutions, we take
one, say $X_3(\zeta)$, with the lowest possible order. And so
forth. At the $N$--step we have a solution $X_N(\zeta)$ of order
$-\kappa_N$. The matrix $X(\zeta)=(X_1(\zeta),X_2(\zeta),\ldots,
X_N(\zeta))$ is called the fundamental matrix of the problem
(\ref{v2}). Vekua \cite{vek} shown that the fundamental matrix
$X(\zeta)$ has the following properties. The determinant of
$X(\zeta)$ is not equal to zero for all $\zeta \in \mathbb{C}$.
Among all fundamental matrices, it is possible to pick one for
which $\lim_{\zeta \to \infty} \zeta^{\kappa_j}X_{jm}(\zeta) =
\delta_{jm}$, where $\delta_{jm}$ is the Kronecker symbol,
$X_{jm}(\zeta)$ is the $m$-th coordinate of the vector function
$X_j(\zeta)$.  The numbers $\kappa_1, \kappa_2,\ldots, \kappa_N$
are called the partial indices. By construction they are ordered
as follows $\kappa_1 \geq \kappa_2 \geq \ldots \geq \kappa_N$. The
sum $\kappa=\sum_{j=1}^N \kappa_j$ is called the total index
(total winding number).

Each solution of (\ref{v2}) in the class
$\mathfrak{H}({\mathbb{D}}^+\cup \dot{\mathbb{D}}^-,V)$ is
represented in the form
\begin{equation}
\label{v3} \Phi(\zeta)=X(\zeta) p(\zeta), \; \zeta \in
(\overline{{\mathbb D}^+} \cup \overline{{\mathbb D}^-})\setminus
V,
\end{equation}
where $p(\zeta)$ is a polynomial vector function. The partial
indices are uniquely determined. In particular, they are
independent of the choice of $X(\zeta)$. These and other
fundamental properties of the fundamental matrix and of the
partial indices are established in \cite{vek}. We describe one of
the most important Vekua's result in the following proposition.

\begin{theorem}[Vekua \cite{vek}, p. 121] \label{ThV}
Let $X(\zeta)$ be a fundamental matrix of the problem (\ref{v1})
and $\kappa_j$ ($j=1,2,\ldots,N$) be its partial indices. Then, the
problem (\ref{v1}) is solvable in the class $\mathfrak{H}_0
(\mathbb{D}^+,V)$ iff
\begin{equation}
\label{4.34}  \int_{|\tau|=1} q(\tau)[X^+(\tau)]^{-1} b(\tau)
d\tau =0,
\end{equation}
where
\begin{equation}
\label{4.q}
q(\tau)=(q_{-\kappa_1-1}(\tau),q_{-\kappa_2-1}(\tau),\ldots,q_{-\kappa_N-1}(\tau)),
\end{equation}
$q_{\nu}(\tau)$ is arbitrary polynomial of degree $\nu$
($q_\nu(\zeta)\equiv 0$ if $\nu<0$).

If (\ref{4.34}) is fulfilled, the general solution of the problem
(\ref{v1}) has the form
\begin{equation}
\label{4.35} \Phi(\zeta) = \frac1{2\pi i} X(\zeta) \int_{|\tau|=1}
\frac{[X^+(\tau)]^{-1} b(\tau)}{\tau - \zeta} d\tau +
X(\zeta)p(\zeta),
\end{equation}
where
\begin{equation}
\label{4.p} p(\zeta)=
(p_{\kappa_1-1}(\zeta),p_{\kappa_2-1}(\zeta),\ldots,p_{\kappa_N-1}(\zeta)).
\end{equation}
Here $p_{\nu}(\zeta)$ is an arbitrary polynomial of degree $\nu \geq 0$, and $p_{\nu}(\zeta) \equiv 0$ for $\nu<0$.
\end{theorem}

In this theorem, Vekua's notations are applied. For instance, the
scalar product of vectors is used in the integrand from
(\ref{4.34}).

\section{Reduction of the Riemann problem to functional equations}
\subsection{Reduction of the Riemann problem to the
Riemann--Hilbert problem}\label{sec2.1} In the present section we
reduce the ${\mathbb C}$--linear conjugation problem (\ref{2.1})
to a boundary value problem. Introduce the vector function
$\phi(z)$ of dimension $2N$
\begin{equation}
\label{2.2} \phi(\zeta) := \left(
\begin{array}{ll}
\Phi(\zeta) \\
\overline{\Phi}\left(\frac1{\zeta}\right)
\end{array}
\right), \quad |\zeta| \leq 1,
\end{equation}
where
$\overline{\Phi}\left(\frac1{\zeta}\right)=\overline{\Phi\left(\frac1{\overline{\zeta}}\right)}$.
The vector function $\phi(\zeta)$ is analytic in $0<|\zeta|<1$,
continuous in $0<|\zeta|\leq 1$ except the points of $V$ where it
is bounded. Moreover, the components $\phi_{N+1}(\zeta),
\phi_{N+2}(\zeta),\ldots , \phi_{2N}(\zeta)$ of the vector
$\phi(\zeta)$ admit a pole of order $M$ at the point $\zeta=0$,
the components $\phi_{1}(\zeta), \phi_{2}(\zeta),\ldots ,
\phi_{N}(\zeta)$ are analytic at $\zeta = 0$. The relation
(\ref{2.1}) can be written in the form
\begin{equation}
\label{2.3} \phi(\zeta)=g_k \overline{\phi(\zeta)}, \; \zeta\in
L_k, \; k=1,2,\ldots,n,
\end{equation}
where
\begin{equation}
\label{2.4} g_k = \left(
\begin{array}{ll}
0 & G_k \\
\overline{G_k^{-1}} & 0
\end{array}
\right).
\end{equation}
By substitution of (\ref{2.2}) and (\ref{2.4}) in (\ref{2.3}) one
can easily see that (\ref{2.3}) is equivalent to (\ref{2.1}).
Hence, the relation (\ref{2.2}) defines an isomorphism between the
problems (\ref{2.1}) and (\ref{2.3}) in the described classes. One
can consider the solution of the problem (\ref{2.1}) as a pair of
the $N$--dimensional vector functions which constitute by
(\ref{2.2}) a
$2N$--dimensional vector function. 

\begin{figure}[hbtp]\centering
\includegraphics[width=0.7\textwidth]{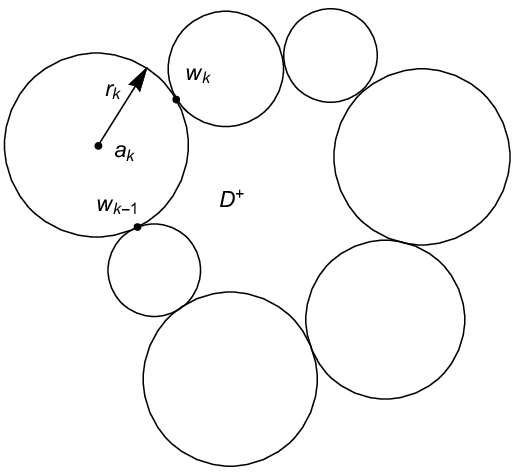}
\caption{Plane $z$. }
\label{geometry}
\end{figure}

Let $n > 2$. Consider a bounded circular polygon $D^+$ with zeroth
angles bounded by touching positively oriented circles $l_k=\{t \in
\mathbb{C}:|t-a_k|=r_k\}$ ($k=1,2,\ldots,n$). We assume that all
points of touching $w_k$ belong to the unit circle. Let $D_k = \{z
\in \mathbb{C}: |z-a_k|<r_k \}$, $D^-$ be the complement of
$(\cup_{k=1}^n \overline{D_k})\cup D^+$ to the extended complex
plane. Divide each circle $l_k$ onto two arcs $l_k^+ = \partial D^+ \cap l_k$ and $l_k^-= \partial D^- \cap l_k$. Let $f:{\mathbb D}^+ \to D^+$ be a conformal mapping of ${\mathbb D}^+$ onto $D^+$. The parameters $a_k$ and $r_k$ are
chosen in such a way that $l_k^+=f(L_k)$. In particular, the
discontinuity set $V$ of the problem (\ref{2.3}) is transformed
onto the touching points of the circles $l_k$, the set
$W=f(V)=\{w_1,w_2,\ldots , w_n\}$, where $w_k = f(\zeta_k)$. The
conformal mapping
\begin{equation}
\label{2.conf}z = f(\zeta)
\end{equation}
can be presented in the form of the Schwarz-Christoffel integral
\cite{neh}. Without loss of generality one can assume that $0 \in D^+$ and take $f(0)=0$. The
symmetry principle implies that
\begin{equation}
\label{2.cob}z =
\left[\overline{f}\left(\frac1{\zeta}\right)\right]^{-1}
\end{equation}
is the conformal mapping of ${\mathbb D}^-$ onto $D^-$ which maps
$\zeta_k$ onto $w_k$. 
%

\begin{remark}
In the case $n=1$, the domain $D^+$ coincides with the unit disk,
i.e., $f(\zeta) = \zeta$. In the case $n=2$, $D^+$ is the strip
$-1 < \Im z < 1$, i.e., $f(\zeta)$ is a composition of the
logarithm and M\"{o}bius transformations. The study of these cases
is much easier than the general case $n>2$. Below, we consider only
the case $n>2$.
\end{remark}

Applying the conformal mapping $z = f(\zeta)$ to (\ref{2.3}) we
arrive at the Riemann--Hilbert problem
\begin{equation}
\label{2.5} \varphi(t)=g_k \overline{\varphi(t)}, \; t\in l_k^+,
\; k=1,2,\ldots,n,
\end{equation}
with respect to the vector function $\varphi=\phi \circ f^{-1}$
analytic in $D^+$, continuous in $\overline{D^+}\setminus W$ and
bounded at the points of $W$. Moreover, the components
$\varphi_{N+1}(z), \varphi_{N+1}(z),\ldots , \varphi_{2N}(z)$ of
the vector $\varphi(z)$ admit a pole of order $M$ at the point
$z=0$.

The inversion $z \mapsto \frac1{\overline{z}}$ transforms the arcs
$l_k^+$ onto $l_k^-$. Consider the auxiliary boundary value
problem
\begin{equation}
\label{2.6} \varphi(t)=g_k \overline{\varphi(t)}, \; t\in l_k^-,
\; k=1,2,\ldots,n,
\end{equation}
with respect to the vector function $\varphi \in
\mathfrak{H}(D^-,W)$. The problems (\ref{2.5}) and (\ref{2.6}) are
related in the following way. Let
$\psi(z)=\overline{\varphi}\left(\frac1z\right)$, $z \in D^+$.
Then, it follows from (\ref{2.6}) that
\begin{equation}
\label{2.7}
\psi(t)=g_k^{-1} \overline{\psi(t)}, \; t\in l_k^+, \;
k=1,2,\ldots,n.
\end{equation}
In accordance with \cite{vek} the problem (\ref{2.7}) is called
the accompanying problem to (\ref{2.5}). Since $g_k \overline{g_k}=I$, $g_k^{-1}$  in \eqref{2.7} can be replaced with $\overline{g_k}$.

\subsection{Reduction of the Riemann--Hilbert problem to the R--linear problem}
Consider the relations (\ref{2.5}) and (\ref{2.6}) as one boundary
value problem with respect to $\varphi \in \mathfrak{H}(D^+ \cup
D^-,W)$ except the components $\varphi_{N+1}(z),
\varphi_{N+1}(z),\ldots , \varphi_{2N}(z)$ which admit a pole of
order $M$ at the point $z=0$.

Consider the matrix $g_k$ defined by (\ref{2.4}). Introduce the
matrix (compare to \cite{MitAnn})
$$
\lambda_k = \left(
\begin{array}{ll}
I & -G_k \\
-iI & -iG_k
\end{array}
\right),
$$
where $I$ is the identity matrix. The matrix $\lambda_k$ is
invertible, since $\det \lambda_k = -2i \det G_k \neq 0$. The inverse
matrix has the form
$$
\lambda_k^{-1} = \frac12 \left(
\begin{array}{ll}
I & iI \\
-G_k^{-1} & iG_k^{-1}
\end{array}
\right).
$$
By direct calculations one can check the following equality
\begin{equation}
\label{3.1} \overline{\lambda_k}^{-1} \lambda_k =-g_k.
\end{equation}
Using (\ref{3.1}) we write (\ref{2.5}) and (\ref{2.6}) in the form
\begin{equation}
\label{3.2} \overline{\lambda_k} \varphi(t)+ \lambda_k \overline{
\varphi(t)}=0, \; t \in l_k, \; k=1,2,\ldots,n,
\end{equation}
where $l_k=l_k^+ \cup l_k^- \cup \{w_{k-1}\}\cup \{w_k\}$.

The boundary value problem (\ref{3.2}) can be written as the
${\mathbb R}$--linear conjugation problem
\begin{equation}
\label{3.3} \overline{\lambda_k} \varphi(t)=\overline{\lambda_k}
\varphi_k(t)- \lambda_k \overline{ \varphi_k(t)}, \quad t \in l_k, \;
k=1,2,\ldots,n,
\end{equation}
where $\varphi_k \in \mathfrak{G}(D_k,\{w_{k-1},w_k\})$.

We now demonstrate that the problems (\ref{3.2}) and (\ref{3.3})
are equivalent. If $\varphi(z)$ and $\varphi_k(z)$ satisfy
(\ref{3.3}), then $\varphi(z)$ satisfies (\ref{3.2}). It is easily
seen by taking the real part from (\ref{3.3}). Conversely, let
$\varphi(z)$ satisfies (\ref{3.2}). Then, we can construct
$\varphi_k(z)$ solving the simple Riemann-Hilbert problem for the
disk $D_k$
\begin{equation}
\label{3.4}
2\Im\;\overline{\lambda_k} \varphi_k(t) = \Im
\overline{\lambda_k} \varphi(t), \quad t \in l_k.
\end{equation}
One can see that (\ref{3.4}) is decomposed onto $2N$ scalar
problems with respect to the components of the vector function
$\overline{\lambda_k} \varphi_k(z)$ and hence are written
explicitly through $\varphi(z)$.
Then $\varphi_k(z)$ is uniquely determined by $\varphi(z)$ up to
an additive constant vector (for details see \cite{gakh}).
Therefore, we have proved that the problems (\ref{3.2}) and
(\ref{3.3}) are equivalent. It is worth noting that $\varphi_k(z)$ are almost bounded at $w_k$ and $w_{k-1}$, e.g., they can have logarithmic singularities at these points. This follows from solution to the problem \eqref{3.4} with respect to $\varphi_k(z)$. The right hand part of \eqref{3.4} can have finite jumps at $w_k$ and $w_{k-1}$. Hence, $\varphi_k(z)$ are almost bounded at $w_k$ and $w_{k-1}$.

Multiplying (\ref{3.3}) by $\overline{\lambda_k}^{-1}$ and using
(\ref{3.1}) we obtain the following ${\mathbb R}$--linear problem
\begin{equation}
\label{3.6}  \varphi(t)= \varphi_k(t)+g_k \overline{
\varphi_k(t)}, \; t \in l_k, \; k=1,2,\ldots,n.
\end{equation}

\subsection{Reduction of the R--linear problem to functional
equations} \label{sec2.3} We now proceed to reduce the problem
(\ref{3.6}) to a system of functional equations.

Introduce the inversion with respect to the circle $l_k$
\begin{equation}
\label{eq:inv}
z_{(k)}^{\ast}=\frac{r_k^2}{\overline{z-a_k}}+a_k.
\end{equation}
It is known that if $\varphi_k(z)$ is analytic in $|z-a_k|<r_k$,
then $\overline{\varphi_k(z_{(k)}^{\ast})}$ is analytic in
$|z-a_k|>r_k$. The functions $\varphi_k(z)$ and
$\overline{\varphi_k(z_{(k)}^*)}$ have the same boundary behavior
on $l_k$ (continuity, boundness), since $t_{(k)}^{\ast}=t$ on
$l_k$.

The orientation of each $l_k$ is counter clockwise. The path along the curve $l$  is defined as follows. It is assumed that the curve $l$ is closed and is consists of the sequent arcs $l_1^+,l_2^+, \ldots, l_n^+$, $l_1^-,l_2^-, \ldots, l_n^-$. Then, the curve $l$ divides the complex plane onto two domains $\mathcal D^+ =\cup_{k=1}^n D_k$ and $\mathcal D^- =D^+ \cup D^-$.

We now proceed to use Sochocki's formulae\footnote{also known as Sokhotsky's and Plemelj's formulae, see for a historical review \cite{mus}} for the closed curve $l$ separating the domains $\mathcal D^{\pm}$. Let $f(t)$ be a scalar function H\"{o}lder continuous separately in each $l_k$ ($k=1,2,\ldots,n$) but having in general a jump at the touching points $w_k$ when $t$ passes from $l_k$ to $l_{k+1}$.
For fixed $m$, introduce the operators
\begin{equation}
\label{eq:la1} S_m f(z) = \frac{1}{2\pi i}\int_{l_m} \frac{f(t)}{t-z} \; dt,  \quad z \in \cup_{k=1}^n D_k \quad (m=1,2,\ldots,n) \footnote{$z$ can belong to $D^+ \cup D^-$}.
\end{equation}
The function $S f(z)=  \sum_{m=1}^n S_m f(z)$ is analytic in the domain $\mathcal D^+ =\cup_{k=1}^n D_k$, continuous in its closure except the touching points where it is almost bounded, more precisely, can have logarithmic singularities \cite{gakh}, \cite{mus}. The function $S f(z)$ satisfies
Sochocki's formulae
\begin{equation}
\label{eq:la2} S f(t) =\lim_{z \to t}S f(z)=\frac{1}{2} \; f(t) + \sum_{m=1}^n\frac{1}{2\pi i}\int_{l_m} \frac{f(\tau)}{\tau-t} \; d\tau,  \quad t \in \cup_{m=1}^n l_m \backslash W.
\end{equation}

Introduce the Banach space $\mathfrak{B}(l,W)$ of vector function
$h(t)=(h_1(t), h_2(t), \ldots, h_{2N}(t))$ H\"{o}lder continuous on each $l_k$ with  the exponent $0<\mu \leq 1$ endowed with the norm
$$
\|h\|=\left( \sum_{j=1}^{2N} \|h_j\|^2 \right)^{1/2}, \quad \|h_j\|=\max_{t\in l} |h(t)|+  \max_{1\leq k \leq n} \sup_{t_1,t_2 \in l_k; t_1\neq t_2} \frac{|h(t_1)-h(t_2)|}{|t_1-t_2|^{\mu}}.
$$
Here, it is assumed that the limit values of $h(t)$ at $w_k \in W$
are the same if $t \to w_k$ with $t \in l_k$, but they can be
different if $t \to w_k$ with $t \in l_k \cup l_{k+1}$. Thus, the
space $\mathfrak{B}(l,W)$ can be also considered as the direct sum
of the spaces of H\"{o}lder continuous vector functions. Vector
functions $h \in \mathfrak{B}(l,W)$ analytically continued into
all disks $D_k$ generate a closed subspace of $\mathfrak{B}(l,W)$
which we denote by $\mathfrak{H}(\cup_{k=1}^n D_k,W)$. The operator $S$ is bounded in $\mathfrak{H}(\cup_{k=1}^n D_k,W)$ \cite{GK}. Instead of $\mathfrak{B}(l,W)$ one can consider the operators in the space $L_2(l)$ \cite{GK}.

The operator $S$ transforms a function $f(t)$ into a function analytic in $\cup_{k=1}^n D_k$.
Apply the operator $S$ to the conjugation condition \eqref{3.6}. First, we have
\begin{equation}
\label{eq:la3f}  S \varphi(z)=\omega(z), \quad z \in D_k, \; k=1,2,\ldots,n,
\end{equation}
where
\begin{equation}
\label{3.6a}  \omega(z)= \left(
\begin{array}{ll}
\gamma_{-1} \\
\overline{\gamma_0}
\end{array}
\right) + \sum_{s=1}^M \left(
\begin{array}{ll}
0 \\
\overline{\gamma_s}
\end{array}
\right) z^{-s}.
\end{equation}
Here, $\gamma_s$ is a constant vector of dimension $N$. The
numeration of $\gamma_s$ and the complex conjugation
$\overline{\gamma_s}$ is taking for unification of the formulae given below.
The vector-function $\omega(z)$ is calculated by residues of $\varphi(z)$ at $z=0$. More precisely, the components $\omega_{N+1}(z), \omega_{N+1}(z),\ldots , \omega_{2N}(z)$ of the vector function $\omega(z)$ can have poles of order $M$ at the point $z=0$, since $\omega(z)$ has the same behavior as $\varphi(z)$ at zero.
We now fix the vector function $\omega(z)$. Equations
on the undetermined vectors $\gamma_s$ ($s= -1,0,\ldots,M$) will
be written latter.

Using the properties of the Cauchy integral we obtain $S \varphi_k(z) = \varphi_k(z)$, $z \in D_k$, and $S \varphi_m(z) = 0$, $z \in D_k$, for $m \neq k$. Therefore, applying the operator $S$ to \eqref{3.6} we arrive at the system of integral equations
\begin{equation}
\label{eq:la4} \varphi_k(z) = \sum_{m=1}^n g_m \;\frac{1}{2\pi i}\int_{l_m} \frac{\overline{\varphi_m(t)}}{t-z} \; dt+\omega(z),  \quad z \in D_k \; (k=1,2,\ldots,n),
\end{equation}
where the matrix $g_m$ is multiplied by the corresponding vector expressed in terms of the Cauchy integral.

The integrals in \eqref{eq:la4} can be calculated by residues. The vector-function  $\overline{\varphi_m(t)}=\overline{\varphi_m(t_{(m)}^*)}$ (see \eqref{eq:inv}) can be analytically continued into $|z-a_m|>r_m$. Hence,
\begin{equation}
\label{eq:la5} \frac{1}{2\pi i}\int_{l_m} \frac{\overline{\varphi_m(t)}}{t-z} \; dt=
\left\{
\begin{array}{ll}
 \overline{\varphi_k(z_{(k)}^*)}- \overline{\varphi_k(a_k)}, \; \mbox{for}\; m=k,
\\
0, \qquad \; \mbox{for} \; m \neq k.
\end{array}
\right.
\end{equation}
Substitution of \eqref{eq:la5} into \eqref{eq:la4} yields the
following system of functional equations with respect to
$\varphi_k(z)$
\begin{equation}
\label{3.7} \varphi_k(z) =  \sum\limits_{m\not=k} g_m
\overline{\varphi_m(z_{(m)}^{\ast})} +\omega(z), \quad z\in
\overline{D_k} \setminus \{w_{k-1},w_k\} \quad  (k=1,2,\ldots,n).
\end{equation}
Here, the constant vectors $\overline{\varphi_k(a_k)}$ are included for shortness in the undetermined vectors $\gamma_{-1}$ and $\gamma_{0}$ of $\omega(z)$.
Let $\varphi_k(z)$ be a solution of (\ref{3.7}). Then the vector
function $\varphi(z)$ is given by the formula
\begin{equation}
\label{3.8}  \varphi(z) =  \sum\limits_{m=1}^{n} g_m
\overline{\varphi_{m}(z_{(m)}^{\ast})}+\omega(z), \quad  z\in
D^+ \cup D^-.
\end{equation}

The special structure \eqref{2.4} of the matrix $g_k$ yields the decomposition of the system \eqref{3.7}
\begin{equation}
\label{eq:3.7a} \Psi_k(z) =  \sum\limits_{m\not=k} G_m
\overline{\widetilde{\Psi}_m(z_{(m)}^{\ast})} +\gamma_{-1},
\end{equation}
\begin{equation}
\label{eq:3.7b} \widetilde{\Psi}_k(z) =  \sum\limits_{m\not=k} \,\overline{G_m^{-1}}
\overline{\Psi_m(z_{(m)}^{\ast})} +\Omega(z), \quad z\in
\overline{D_k} \setminus w_k,w_{k+1} \quad  (k=1,2,\ldots,n),
\end{equation}
where the $2N$-dimensional vector functions $\varphi_k(z)$ are presented through two $N$-dimensional vector functions $\Psi_k(z)$ and $\widetilde{\Psi}_k(z)$
\begin{equation}
\label{eq:3.7c}
\varphi_k(z) = \left(
\begin{array}{ll}
   \Psi_k(z) \\
\widetilde{\Psi}_k(z)
\end{array}
\right)
\end{equation}
and
\begin{equation}
\label{eq:3.7d}
\Omega(z)=\sum_{s=0}^M \overline{\gamma_s} z^{-s}.
\end{equation}
The vector functions $\widetilde{\Psi}_k(z)$ can be eliminated from equations \eqref{eq:3.7a}-\eqref{eq:3.7b}
\begin{equation}
\label{eq:3.7e}
 \Psi_k(z) =  \sum\limits_{m\not=k}  \sum\limits_{\ell\not=m} G_m G_{\ell}^{-1}
\Psi_{\ell}(z_{(\ell m)}^{\ast}) +F_k(z), \quad z\in
\overline{D_k} \setminus w_k,w_{k+1} \quad  (k=1,2,\ldots,n),
\end{equation}
where
$$
z_{(\ell m)}^{\ast} =\left(z_{(m)}^{\ast} \right)_{(\ell)}^{\ast} =
\frac{r_l^2(z-a_m)}{r_m^2-(\overline{a_l-a_m})(z-a_m)}+a_l
$$
denotes the compositions of inversions with respect to the $m$th and the $\ell$th circles and
\begin{equation}
\label{eq:3.7f}
F_k(z)=\gamma_{-1}+\sum\limits_{m\not=k} G_m \overline{\Omega\left(z_{(m)}^{\ast} \right)}=
\gamma_{-1}+\sum_{s=0}^M \sum\limits_{m\not=k} G_m \gamma_s \left(\overline{z_{(m)}^{\ast}} \right)^{-s}.
\end{equation}
One can see that each $F_k(z)$ is the linear combination of the  functions $\{1, \left(\overline{z_{(m)}^{\ast}} \right)^{-s} \}$ with the coefficients consisting of the components of $\gamma_s$. Hence, $F_k(z)$ is analytic in $\overline{D_k}$. 

Write equations \eqref{eq:3.7e} in the form
\begin{equation}
\label{eq:3.7g}
\begin{array}{ll}
 \Psi_k(z) =  G_{k-1} G_{k}^{-1} \Psi_{k}(z_{(k,k-1)}^{\ast})+G_{k+1} G_{k}^{-1} \Psi_{k}(z_{(k,k+1)}^{\ast})
+\sum\limits_{(m, \ell) \in \mathcal N_k} G_m G_{\ell}^{-1}
\Psi_{\ell}(z_{(\ell m)}^{\ast}) +F_k(z), \\ \qquad \qquad \qquad \qquad z\in
\overline{D_k} \setminus w_k,w_{k+1}\qquad  (k=1,2,\ldots,n),
\end{array}
\end{equation}
where $k-1=n$ for $k=1$.
The set $\mathcal N_k$ consists of the pairs $(m, \ell)$ with not equal elements $m, \ell \in \mathbb Z_n$; the pairs with $m=k$ and two pairs $(k,k-1)$, $(k,k+1)$ are also eliminated.
The shifts $z_{(\ell m)}^{\ast}$ are M\"{o}bius transformations of elliptic type \cite{Ford}. For $(m, \ell) \in \mathcal N_k$, every $z_{(\ell m)}^{\ast}$  transforms the disks $|z-a_k| \leq r_k$ and $|z-a_\ell| \leq r_\ell$ into the disk $|z-a_\ell| < r_\ell$. The map $z_{(\ell m)}^{\ast}$ has a unique attractive fixed point $z_*$ in $|z-a_\ell| < r_\ell$ \cite{Ford, MiR}.

 The M\"{o}bius transformation $z_{(k,k-1)}^{\ast}$ and $z_{(k,k+1)}^{\ast}$ are of parabolic type \cite{Ford} with the neutral fixed points $w_{k-1}$ and $w_k$, respectively, since the $(k-1)$th, $k$th circles touch at $w_{k-1}$ and $(k+1)$th, $k$th at $w_{k+1}$.

Equations \eqref{eq:3.7g} can be considered in the Hardy type space $\mathcal{H}_2(\cup_{k=1}^n D_k)$ associated with $L_2$. Let $\Psi_k(z)$ belongs to the Hardy space $\mathcal{H}_2(D_k)$, i.e., $\Psi_k(z)$ is analytic in $D_k$ and the norm is defined by formula
\begin{equation}
\label{eq:har1} \|\Psi_k\|_2 = \sup_{0<r<r_k} \left[ \frac{1}{2\pi} \int_0^{2\pi} |\Psi_k(r e^{i \theta})|^2 d \theta \right]^{\frac 12}.
\end{equation}
A function $\Psi(z)$ belongs to  $\mathcal{H}_2(\cup_{k=1}^n D_k)$ if $\Psi \in \mathcal{H}_2(D_k)$ for all $k=1,2, \ldots, n$ and $\|\Psi\|_2 = \max\limits_{1\leq k \leq n} \|\Psi_k\|_2$. The norm of the vector functions is introduced by the norms of components as $\|\Psi\|_2=\left( \sum_{j=1}^{N} \|\Psi_j\|_2^2 \right)^{1/2}$.
It was shown in \cite{MiR} that the operator $\sum\limits_{(m, \ell) \in \mathcal N_k} G_m G_{\ell}^{-1} \Psi_{\ell}(z_{(\ell m)}^{\ast})$ is compact in the spaces $\mathfrak{H}(\cup_{k=1}^n D_k,W)$ and $\mathcal{H}_2(\cup_{k=1}^n D_k)$.

\begin{theorem}
Let $\Psi \in \mathcal{H}_2(\cup_{k=1}^n D_k)$ where $\Psi(z)=\Psi_{k}(z)$ in $D_k$, i.e., the vector function $\Psi(z)$ is separately defined in all the disks $D_k$.
The operators
$\mathcal{A} \Psi(z)= G_{k-1} G_{k}^{-1} \Psi_{k}(z_{(k,k-1)}^{\ast})$ and $\mathcal{B} \Psi(z)=G_{k+1} G_{k}^{-1} \Psi_{k}(z_{(k,k+1)}^{\ast})$, $z\in D_k$ ($k=1,2,\ldots,n$) are compact in $\mathcal{H}_2(\cup_{k=1}^n D_k)$.
\end{theorem}

Proof.  The operators $\mathcal{A}$ and $\mathcal{B}$ consist of the compositions of the multiplications by matrices and of the shift operators. The multiplication operators are bounded.
It is sufficient to investigate the scalar shift operators. For definiteness, consider the operator
\begin{equation}
\label{eq:a8}
\mathcal{F}f(z)= f(z_{(k,k-1)}^{\ast}).
\end{equation}
The general form of the map $z_{(k,k-1)}^{\ast}$ can be reduced by translations and rotations to the case $w_k=0$, $a_k=r_1>0$, $a_{k-1}=-r_2<0$. The compactness property of $\mathcal{F}$ does not change after such transformations. In this case,
\begin{equation}
\label{eq:4.1a}
z_{(k,k-1)}^{\ast} = \frac{z}{1+\left(\frac{1}{r_1}+\frac{1}{r_2} \right) z}=:\alpha(z).
\end{equation}
Compactness of the operator \eqref{eq:a8} follows from the Hilbert-Schmidt Theorem for composition operators \cite[page 26]{Shapiro}. It is sufficient to check the boundness of the integral
\begin{equation}
\label{eq:4.1b}
\int_{-\pi}^{\pi} \frac{d\theta}{1-|\frac{1}{r_1} \alpha[r_1(1+e^{i \theta})]|^2} =2 \pi \;\frac{r^2+3r+2+ \frac{1}{\sqrt{4 r^2+12   r+13}}}{r^2+3r+3},
\end{equation}
where $r=\frac{r_1}{r_2}$. The integral \eqref{eq:4.1b} is computed with Mathematica$^\circledR$.

The theorem is proved.

\section{Solution to functional equations and calculation of partial indices}
\label{pi}
The space $\mathcal{H}_2(\cup_{k=1}^n D_k)$ of scalar functions is isomorphic to the space $h_2$ of the sequences $\{\{\psi_{kj} \}_{k=1}^n\}_{j=0}^{\infty}$ with the norm 
$$
\|\{\psi_{kj} \}\|_2=\sum\limits_{k=1}^{n} \left( \sum\limits_{j=0}^{\infty}|\psi_{kj}|^2 \right)^{1/2}
$$ 
where the Taylor expansion of  $\Psi_k(z)$ is used
\begin{equation}
\label{eq:har2}
\Psi_k(z) =\sum_{j=0}^{\infty} \psi_{kj} (z-a_k)^j, \quad |z-a_k|<r_k.
\end{equation}
Here, $\psi_{kj}=\frac{1}{j!} \Psi^{(j)}_k \Big|_{z=a_k}$ denotes the derivative of order $j$ of the function $\Psi_k$ at $z=a_k$.
The spaces of vector functions are introduced by coordinates. It is worth noting that $\mathfrak{H}(\cup_{k=1}^n D_k,W) \subset \mathcal{H}_2(\cup_{k=1}^n D_k)$ and given by \eqref{eq:3.7f} the terms $F_k(z)$ in equations \eqref{eq:3.7g} belong to $\mathfrak{H}(\cup_{k=1}^n D_k,W)$. It will be shown later that any solution in $\mathcal{H}_2(\cup_{k=1}^n D_k)$ with such $F_k(z)$ belongs to $\mathfrak{H}(\cup_{k=1}^n D_k,W)$, e.g. Pumping Principle \cite[page 22]{MiR}. Therefore, one can formally solve the problem in $\mathcal{H}_2(\cup_{k=1}^n D_k)$. The obtained solution will belong to the space $\mathfrak{H}(\cup_{k=1}^n D_k,W)$.

Using the general properties of equations with compact operators in the Hilbert space $h_2$ \cite{AG,KK} we can develop the following constructive algorithm to solve the system \eqref{eq:3.7g} and to calculate the partial indices of the corresponding Riemann-Hilbert problem.  We have
\begin{equation}
\label{eq:har3}
\Psi_{\ell}(z_{(\ell m)}^{\ast}) =\sum_{j=0}^{\infty} \psi_{lj} \left[\frac{r_l^2(z-a_m)}{r_m^2-(\overline{a_l-a_m})(z-a_m)}\right]^j, \quad |z-a_k|<r_k,
\end{equation}
The discrete form of the functional equations is obtained after substitution of \eqref{eq:har2}-\eqref{eq:har3} into \eqref{eq:3.7e}
\begin{equation}
\label{eq:syst}
\begin{array}{ll}
 \sum\limits_{j=0}^{\infty} \psi_{kj} (z-a_k)^j = \sum\limits_{m\not=k}  \sum\limits_{\ell\not=m} G_m G_{l}^{-1}\sum\limits_{j=0}^{\infty} \psi_{lj} \left[\frac{r_l^2(z-a_m)}{r_m^2-(\overline{a_l-a_m})(z-a_m)}\right]^j+\gamma_{-1}
\\
+\sum\limits_{s=0}^M \sum\limits_{m\not=k}^n G_m \gamma_s \left(\frac{r_m^2}{z-a_m}+\overline{a_m} \right)^{-s},  \quad z\in
\overline{D_k} \setminus \{w_k,w_{k+1}\} \quad  (k=1,2,\ldots,n).
\end{array}
\end{equation}
All the series in \eqref{eq:syst} can be expand in $(z-a_k)^j$ in the disk $|z-a_k|<r_k$.
Comparison of the coefficients on $(z-a_k)^j$ yields an infinite system of linear algebraic equations on the coefficients $\psi_{lj}$ and parameters of $\omega(z)$. The infinite system can be solved by the truncation method \cite{KK}.
Taking a partial sum of Taylor series \eqref{eq:har2} with $Q$ first items
and collecting the coefficients on the same powers of $z-a_k$
we arrive at the finite system on $\psi_{kj}$
\begin{equation}
\label{eq:syst2}
\psi_{kj}= \sum\limits_{m \neq k}^n\sum\limits_{l \neq m}^n\sum\limits_{i=0}^Q \mu_{ml}^{(kij)} \,G_m G_l^{-1}\psi_{li}  +F_k^{(j)}(a_k), \quad
k=1,\dots,n,\; j=0,\dots,Q,
\end{equation}
where the scalar
\begin{equation}
\label{syst1}
\mu_{ml}^{(kij)}=\frac{1}{j!}\left[\left( \frac{r_l^2(z-a_m)}{r_m^2-(\overline{a_l-a_m})(z-a_m)} \right)^i\right]^{(j)}\Big|_{z=a_k},
\end{equation}
and the vectors
\begin{equation}
\label{eq:syst3}
F_k^{(j)}(a_k)=\frac{1}{j!}\left[\gamma_{-1}+\sum\limits_{s=0}^M \sum\limits_{m\not=k}^n G_m \gamma_s \left(\frac{r_m^2}{z-a_m}+\overline{a_m} \right)^{-s}\right]^{(j)}\Big|_{z=a_k}.
\end{equation}

The vectors $\gamma_s$ can be also considered as unknowns. Then, \eqref{eq:syst2} can be considered as a homogeneous system of equations on $\boldsymbol{\Psi}=\{\psi_{kj}\}_{k=1,\dots,n;j=0,\dots}$ and $\boldsymbol{\gamma}=\{\gamma_{-1}, \gamma_0,\cdots \gamma_M \}$.
The left hand part of \eqref{eq:syst2} determines a compact operator $\mathbf A$ in the space $h_2$. Then, the system \eqref{eq:syst2} can be shortly written as the operator equation in $h_2$
\begin{eqnarray}
\label{eq:S1}
\boldsymbol{\Psi} = \mathbf A \boldsymbol{\Psi} + \mathbf F \boldsymbol{\gamma},
\end{eqnarray}
where the finite rank operator $\mathbf F$ is determined by \eqref{eq:syst3}.
The truncated system is denoted as
\begin{eqnarray}
\label{eq:S2}
\boldsymbol{\Psi}_Q = \mathbf A_Q \boldsymbol{\Psi}_Q + \mathbf F \boldsymbol{\gamma},
\end{eqnarray}
where $\mathbf A_Q$ is the truncated matrix of the dimension $n N(Q+1) \times n N(Q+1)$, $\boldsymbol{\Psi}_Q=\{\psi_{kj}\}_{k=1,\dots,n;j=0,\dots,Q}$ is a column vector with $n N(Q+1)$ entries. Solution to \eqref{eq:S2} approximates solution to the infinite system \eqref{eq:S1}  for sufficiently large $Q$ \cite{AG,KK}.
This fact is summarized in the following theorem which relates solution to the Riemann-Hilbert problem with linear systems.

\begin{theorem}
\label{thF} Let $G_k$ be invertible matrices ($k \in
\mathbb{Z}_n$, $n >2$) and each matrix $G_k^{-1}G_{k-1}$ has
different eigenvalues. The problem (\ref{2.1}) in
$\mathfrak{H}({\mathbb{D}}^+\cup \dot{\mathbb{D}}^-,V)$ with the
prescribed order at infinity $M$ is solvable if and only if the
linear algebraic system \eqref{eq:S1} on the
$N$--vectors $\psi_{kj}$ and $\gamma_s$ ($s=-1,0, \ldots, M$) is solvable.

Solvability of the system \eqref{eq:S1} is reduced to solvability of the system \eqref{eq:S2}  for sufficiently large $Q$. Its solution can be approximated by solution to the system \eqref{eq:S2}.  

Let the general solution to the finite system \eqref{eq:S2} depend  on the arbitrary vectors $\gamma_s$ ($s=r,r+1,\ldots, M$), i.e., the rank of the system \eqref{eq:syst2} is equal to $r+1$. Then the problem (\ref{2.1}) has $M-r+1$ linearly
independent solutions approximated by polynomial corresponding to the truncated series \eqref{eq:har2}.
\end{theorem}

Applying directly the definition of partial indices described in Sec.\ref{sec2b} and Theorem \ref{thF} to the problem  (\ref{2.1}) one can calculate the partial
indices. Considering $M$ as a control parameter we investigate the linear system \eqref{eq:S2} for sufficiently large $Q$.

Let $M=0$. If the system \eqref{eq:S1} has non--zero solutions, then $\kappa_1=0$ and the first column $X_1(\zeta)$ of the fundamental matrix $X(\zeta)$ is the constant
vector $\gamma_{-1}$ in $|\zeta| \leq 1$ and $\gamma_{0}$ in
$|\zeta| \geq 1$. Such a case is possible since (\ref{2.1}) is reduced to
\begin{equation}
\label{7.1} \gamma_{-1}=G_k \gamma_0, \quad k=1,2, \ldots, n
\end{equation}
with $\gamma_0 = (1,0,\ldots,0)^\top$. Among all matrices $G_k$ at
least two ones are different, since $n>2$ and $G_k^{-1}G_{k-1}$
has different eigenvalues. Equations (\ref{7.1}) on $\gamma_{-1}$ and $\gamma_0$ can have non--trivial solutions. For instance,
$G_k= \left(
\begin{array}{ll}
1 \quad 0 \\
0 \quad k
\end{array}
\right)$, $k=1,2, \ldots, n$ and $\gamma_{-1}=\gamma_0=\left(
\begin{array}{ll}
1 \\
0
\end{array}\right)$ satisfy (\ref{7.1}). 
Here, one can see the principal difference between the scalar ($N=1$) and vector ($N>1$) cases, because for $n>1$ the scalar problem (\ref{2.1}) always has only zero solution.

If the system \eqref{eq:S2} has only zero solution, we increase $M$ up to 1 and consider the corresponding system with $\gamma_{-1}$, $\gamma_0$ and $\gamma_1$. If this system has non--zero solutions, the partial index $\kappa_1=-1$ and $X_1(\zeta)$ can be approximately found. If the system \eqref{eq:S2} with $M=1$ has only zero solution, we again increase $M$ and so forth. Therefore, the partial indices are precisely calculated by the ranks of the system \eqref{eq:S2} for different $M$. It follows from the general theory this process of increasing $M$ will be always finite.

\section{Discussion}
In the present paper, the vector-matrix  Riemann boundary value
problem (${\mathbb C}$--linear conjugation problem) for the unit
disk with piecewise constant matrix has been solved in H\"{o}lder and Hardy spaces by the method of functional equations.


First, the conformal mapping (\ref{2.conf}) $z = f(\zeta)$ of the
unit disk onto a circular polygon is not constructed here. It is
worth noting that this problem differs from the classic problem of
inverse mapping. Our problem is easier, since the points $\zeta_k$
are known. What is more, we suggest that any Christoffel-Schwarz
transformation of the unit disk onto a circular polygon bounded by
externally touching circles, where all touching points lie on one
circle, can be taken as $f(\zeta)$.


The solution of the problem can be found in the form of the series. The second part will contain numerical examples of the above constructive scheme.

\section*{Acknowledgments}
I am grateful to Prof. Frank-Olme Speck for valuable discussion during my visit in 2006 supported by Centro de Matem\'{a}tica e
Aplica\c{c}\~{o}es of Instituto Superior T\'{e}cnico, Lisboa. I thank Ekaterina Pesetskaya and Gia Giorgadze for fruitful discussions and preliminary numerical examples not presented in this paper.


\begin{thebibliography}{99}

\bibitem{AG}
N.I. Akhiezer and I.M. Glazman, Theory of Linear Operators in Hilbert Space,  New York, Dover Publications, 1993.

\bibitem{Antipov}
Y.A. Antipov, Vector Riemann-Hilbert problem with almost periodic and meromorphic coefficients and applications, Proc. R. Soc. A. (2015) 471, 20150262.

\bibitem{ErhSpe02}
T. Ehrhardt, F.-O. Speck, Transformation techniques towards the
factorization of 2$\times$2  matrix functions, {\it Linear Algebra
Appl.}, {\bf 353}, No. 1-3, 53-90, 2002.

\bibitem{ErhSp01}
T.Ehrhardt, I.M.Spitkovsky, Factorization of piecewise constant
matrix-valued functions, and systems of linear differential
equations, {\it St.Petersbg. Math. J.}, {\bf 13}, No. 6, 939-991,
2001.

\bibitem{Ford} L. Ford, Automorphic Functions, AMS, Rhode Island, 1951.

\bibitem{Gia}
G. Giorgadze, G. Khimshiashvili, The Riemann-Hilbert problem in loop spaces,
Doklady Mathematics 73 (2), 258-260, 2006.

\bibitem{GiaA}
G.K. Giorgadze, Solvability condition of the Riemann-Hilbert monodromy problem, Progress In Analysis And Its Applications: Proceedings of the 7th International ISAAC Congress,  ed. M. Ruzhansky, Singapore,  89-95, 2010.

\bibitem{Gia2}
G. Giorgadze, Some analytical and geometrical aspects of stable partial indices,
J. Math. Sci. 193 (3), 461-477, 2013.

\bibitem{Gia3}
G. Giorgadze, N. Manjavidze, On some constructive methods for the matrix Riemann-Hilbert boundary-value problem, J. Math. Sci. 195 (2), 146-174,2013.

\bibitem{GK}
I. Gohberg, N. Krupnik, One-dimensional Linear Singular Integral Equations, v.1-2, Birkh\"{a}user Basel, 1992.

\bibitem{gakh} F. D. Gakhov, {\it Boundary value problems}, Pergamon Press, Oxford, 1966; Nauka, Moscow, 1977.

\bibitem{KK}
L.V. Kantorovich, V.I. Krylov, Approximate Methods of Higher Analysis,  Groningen 1958.

\bibitem{Khr71}
A.A. Khrapkov, Certain cases of the elastic equilibrium of an
infinite wedge with a non-symmetric notch at the wertex, subject
to concentrated force, {\it PMM., J. Appl. Mat. Mech.}, {\bf 35},
625-637, 1971.

\bibitem{Khv94}
L.A. Khvoshchinskaya, The homogeneous boundary value Riemann
problem for two pairs of functions with the piecewise constant
matrix in the case of four and more singular points, {\it Izv.
Akad. Nauk Belarusi, Ser. Fiz.-Mat. Nauk}, No. 2, 124-125, 1994
(summary of the article deposited at VINITI, No. 2400-B93), in
Russian.

\bibitem{kra} M. A. Krasnoselskij, G. M. Vainikko, P. P. Zabreiko, Ya. B. Rutitskij,
V. Ya. Stetsenko, {\it Approximate Solution of Operator
Equations}, Wolters-Noordhoff: Groningen, 1972.

\bibitem{LitSp87}
G. S. Litvinchuk, I. M. Spitkovsky, {\it Factorization of
Measurable Matrix Functions}, Oper. Theory: Adv. Appl. {\bf 25},
Birkh\"auser, Basel, 1987.

\bibitem{MeiSpe89}
E. Meister, F.-O. Speck, The explicit solution of elastodynamical
diffraction problems by symbol factorization, {\it Z. Anal. Anw.}
{\bf 8}, 307-328, 1989.

\bibitem{mikh}
S. G. Mikhlin, {\it Integral Equations and Their Applications to
Certain Problems in Mechanics, Mathematical Physics and
Technology}, Pergamon, New York, 1957.

\bibitem{MitAnn}  V. Mityushev, Riemann problem on double of multiply
connected region, Annales Polon. Math., 1997, 77.1, 1-14.

\bibitem{MiR} V. V. Mityushev,  S. V. Rogosin, {\it Constructive methods for linear and
non-linear boundary value problems for analytic function. Theory
and applications}, Chapman \& Hall / CRC, Monographs and Surveys
in Pure and Applied Mathematics, Boca Raton etc, 1999.

\bibitem{RogMit} V. V. Mityushev,  S. V. Rogosin, {\it On the Riemann-Hilbert problem with a piecewise constant matrix}, ZAA, 27, 53-66, 2008.

%

\bibitem{mus} N. I. Muskhelishvili, {\it Singular Integral Equations: Boundary
Problems of Function Theory and Their Applications to Mathematical
Physics}, New York: Dover Publications; 2nd ed., 1992.

\bibitem{neh} Z. Nehari, {\it Conformal Map}, New York: Dover, 1982.
Driscoll, Tobin A.; Trefethen, Lloyd N. (2002), Schwarz-Christoffel mapping, Cambridge Monographs on Applied and Computational Mathematics 8, Cambridge University Press.

\bibitem{Shapiro}
J.H. Shapiro, Composition operators and classical function theory,
Springer-Verlag, New York, 1993.

\bibitem{vek}
N. P. Vekua, {\it Systems of Singular Integral Equations}, Groningen: Noordhoff, 1967; Nauka: Moscow, 2nd. ed., 1970 [in Russian].

\bibitem{ZveKhv85}
E. I. Zverovich, L. A. Khvoshchinskaya, Solution of a generalized
Riemann problem, In: {\it Scientific works of the Jubilee seminar
on boundary value problems}, Collect. artic. (V.P.Platonov et al.,
eds.), University Publ., Minsk, 1985, 153-158.

\end{thebibliography}
\end{document}